\documentclass{amsart}

\usepackage{amsmath,amssymb,amsthm,a4wide}
\usepackage{graphicx,tikz}

\newtheorem*{thm}{Theorem}

\newtheorem*{corollary}{Corollary}

\theoremstyle{definition}

\theoremstyle{remark}

\begin{document}

\title[]{A Compactness Principle for Maximizing \\Smooth Functions over Toroidal Geodesics}

\author[]{Stefan Steinerberger}
\address{Department of Mathematics, Yale University, New Haven, CT 06511, USA}
\email{stefan.steinerberger@yale.edu}

\thanks{This work is supported by the NSF (DMS-1763179) and the Alfred P. Sloan Foundation.}

\begin{abstract} Let $f \in C^2(\mathbb{T}^2)$ have mean value 0 and consider 
$$ \sup_{\gamma~{\tiny \mbox{closed geodesic}}}{~~~ \frac{1}{|\gamma|} \left| \int_{\gamma}{ f ~~d\mathcal{H}^1}\right| },$$
where $\gamma$ ranges over all closed geodesics $\gamma:\mathbb{S}^1 \rightarrow \mathbb{T}^2$ and $|\gamma|$ denotes their length. We prove that this supremum is always attained. Moreover, we can bound the length of the geodesic $\gamma$ attaining the supremum in terms of \textit{smoothness} of the function: for all $s \geq 2$,
$$ |\gamma|^{s} \lesssim_s  \left( \max_{|\alpha| = s}{ \| \partial_{\alpha} f \|_{L^{1}(\mathbb{T}^2)}} \right) \| \nabla f \|_{L^2}^{} \|f\|_{L^2}^{-2}.$$
This seems like an interesting phenomenon. We do not know at what level of generality it holds or whether versions of it could be established in other settings (hyperbolic surfaces, groups,...).
\end{abstract}

\maketitle

\section{Introduction and main result}
The purpose of this short note is to discuss an interesting phenomenon: let $f$ be a smooth function with mean value 0 on $\mathbb{T}^2$ and suppose we are interested in finding the largest (absolute) average value $f$ can assume on closed geodesics. 
More precisely, if $\gamma:\mathbb{S}^1 \rightarrow \mathbb{T}^2$ is a closed geodesic on $\mathbb{T}^2$, then we are interested in the maximal possible size of
 $$\frac{1}{|\gamma|} \left| \int_{\gamma}{ f ~~d\mathcal{H}^1}\right|,$$
where $\mathcal{H}^1$ is the Hausdorff measure or, since everything is smooth, the usual arclength measure and $|\gamma|$ is the length of the geodesic.
Our result states that the maximal value is assumed for a geodesic of finite length and \textit{that length can be bounded by the smoothness of the function}. 

\begin{thm} Let $f:\mathbb{T}^2 \rightarrow \mathbb{R}$ be at least $s \geq 2$ times differentiable and have mean value 0. Then
$$ \sup_{\gamma~{ \small \emph{closed geodesic}}}{~~~ \frac{1}{|\gamma|} \left| \int_{\gamma}{ f~~d\mathcal{H}^1 }\right| },$$
is assumed for a closed geodesic $\gamma:\mathbb{S}^1 \rightarrow \mathbb{T}^2$ of length no more than
$$ |\gamma|^{s} \lesssim_s  \left( \max_{|\alpha| = s}{ \| \partial_{\alpha} f \|_{L^{1}(\mathbb{T}^2)}} \right) \| \nabla f \|_{L^2}^{} \|f\|_{L^2}^{-2}.$$
\end{thm}
We believe this to be a rather interesting phenomenon. A priori, it might seem very implausible: clearly we can take a closed geodesic $\gamma$ as long as we like, define $f$ to be, say, 1 on the geodesic 
and have it assume smaller values everywhere else (while balancing it in such a way that the mean value is 0). The main result states that this cannot be done without either introducing very
large derivatives or having the function be so negative so as to create shorter geodesics assuming larger extreme values somewhere else. Or, put differently, smoothness of the function $f$ is enough to ensure simplicity
of the extremizing geodesic. We have no reason to believe that the inequality is sharp but the result is not arbitrarily far away from the truth: if we consider $f(x,y) = \sin{(x + \ell y)}$ for some $\ell \in \mathbb{N}$, then the extremizing geodesic
has length $|\gamma| \sim \ell$ while
$$ \max_{|\alpha| = s}{ \| \partial_{\alpha} f \|_{L^{1}(\mathbb{T}^2)}}\sim \ell^s, \quad \| \nabla f \|_{L^2}^{} \sim \ell \qquad \mbox{and} \qquad \|f\|_{L^2}^{-2} \sim 1.$$
For examples of this type, functions whose Fourier series is compactly supported, we can explicitly control the limit $s \rightarrow \infty$ and obtain
a sharp result.

\begin{corollary} Let $N \in \mathbb{N}$ and let $f:\mathbb{T}^2 \rightarrow \mathbb{R}$  be of the form
$$ f(x) = \sum_{\|k\| \leq N}{\widehat{f}(k) e^{2\pi i \left\langle k, x \right\rangle}}.$$
Then the supremum
$$ \sup_{\gamma~{\tiny \emph{closed geodesic}}}{~~~ \frac{1}{|\gamma|} \left| \int_{\gamma}{ f ~~d\mathcal{H}^1}\right| }$$
is attained by a closed geodesic $\gamma$ with $|\gamma| \lesssim N$.
\end{corollary}

The main purpose
of this paper is to introduce the phenomenon and ask a simple question.

\begin{quote}
\textbf{Question.} Does such a principle exist on more general compact manifolds? Are there examples of other geometries where such a bound can be established?
\end{quote}

We are not aware of any results in this direction. Variants of our results can be established on $\mathbb{T}^d$ with $d \geq 3$ but the mechanism is the same and will not yield additional insight into whether this compactness phenomenon is true in a general context or on other non-toroidal geometries.
A result on very general manifolds may either be false or is likely out of reach since already structural statements about closed geodesics are highly nontrivial (see Berger \cite[\S 10.4]{berger} or Klingenberg \cite{klingenberg}). Possible candidates for examples might be hyperbolic surfaces on which there exists a suitably accessible description of closed geodesics or groups on which Fourier Analysis is well understood.

\section{Proofs}
\subsection{Proof of the Theorem}
\begin{proof} The main idea of the proof is rather simple: we show that there exists a relatively short geodesic for which the absolute value of the arising mean value is at least of a certain size. This is done with an averaging argument and non-constructive. The second part of the proof shows that all long geodesics are uniformly bounded. Throughout the paper, we use $A \lesssim B$ to denote $A \leq cB$ for some universal constant $c>0$ and $A \lesssim_s B$ to denote that the implicit constant depends on $s$.
 We identify $\mathbb{T}^2 \cong [0,1]^2$ and write
$$ f(x) = \sum_{ k \in \mathbb{Z}^2}{ \widehat{f}(k) e^{2 \pi i \left\langle x, k \right\rangle} }.$$
Since we assume $f \in C^s(\mathbb{T}^2)$ with $s \geq 2$, the Fourier series converges. Moreover, since $f$ has mean value 0, we have $\widehat{f}(0,0) = 0$. 
Closed geodesics on $\mathbb{T}^2$ can be written as 
$$\gamma(t) = (at, bt + c),\quad  \mbox{where}~(0,0) \neq (a,b) \in \mathbb{Z}^2 ~\mbox{and}~ 0 < c < 1.$$
We assume henceforth that $\mbox{gcd}(a,b) = 1$, this implies that the closed geodesic can be parametrized by $0 \leq t \leq 1$. Moreover, this induces a bijective mapping between closed geodesics and the set of coprime pairs of integers $(a,b) \times [0,1]$. We now argue that it is possible to assume w.l.o.g. that
$$ \sum_{k \in \mathbb{Z}}{ | \widehat{f}(k,0) |^2} \leq \frac{\|f\|^2_{L^2}}{2}.$$
If this were false, then we consider the function $\tilde f(x,y) = f(y,x)$:  geodesics over $f$ correspond to geodesics over $\tilde f$ and Plancherel's theorem implies that the desired inequality is now satisfied. The integral over a closed geodesic can be written as
\begin{align*}
 \frac{1}{|\gamma|} \int_{\gamma}{f ~~d\mathcal{H}^1} &= \int_{0}^{1}{f(\gamma(t)) dt} = \int_{0}^{1}{ \sum_{k \in \mathbb{Z}^2}{ \widehat{f}(k) e^{2 \pi i \left\langle \gamma(t), k \right\rangle} } }  =  \int_{0}^{1}{ \sum_{k \in \mathbb{Z}^2}{ \widehat{f}(k) e^{2 \pi i (k_1 a t + k_2 (bt +c)) }} dt} \\
&=    \int_{0}^{1}{  \sum_{k \in \mathbb{Z}^2}{ \widehat{f}(k) e^{2\pi i k_2 c}  e^{2 \pi i t (k_1 a + k_2 b) } }  dt}  =  \sum_{k_1 a + k_2 b = 0}{ \widehat{f}(k_1, k_2) e^{2\pi i k_2 c} }.
\end{align*}
We now interpret this as a function in $c$. For any fixed $(0,0) \neq (a,b) \in \mathbb{Z}^2$ and coprime $a,b$, the only solutions to $ k_1 a + k_2 b = 0$ are given by $(k_1,k_2) = (-d b, d a)$ where $d \in \mathbb{Z}$.
We will only consider $(0,0) \neq (a,b)$ for which $a \neq 0$: if $a=0$, then the representation for the integral over a closed geodesic does not give rise to a Fourier series but collapses to a number instead.
This is not surprising considering that $\gamma(t) = (at, bt + c)$: if $a = 0$, then the closed geodesic has length 1 and is invariant under $c$.
An application of Plancherel's theorem yields
$$
 \sum_{k \in \mathbb{Z}^2 \atop k_1 a + k_2 b = 0}{ |\widehat{f}(k_1, k_2)|^2} = \left\|  \sum_{ k_1 a + k_2 b = 0}{ \widehat{f}(k_1, k_2) e^{2\pi i k_2 c} } \right\|^2_{L^2[0,1]} \leq 
\max_{0 \leq c \leq 1}{~~ \left|  \frac{1}{|\gamma_{a,b,c}|} \int_{\gamma_{a,b,c}}{f~~ d\mathcal{H}^1}  \right|^2}.$$
Let now $N \in \mathbb{N}$ be a large number to be fixed later. We will sum this inequality over all geodesics $\gamma(t) = (at, bt + c)$ with $|a|, |b| \leq N$ and $a \neq 0$ in such a way that the arising lines cover all lattice points outside the $x-$axis in a radius $\sim N$. All these geodesics have length $\sqrt{a^2+b^2} \leq 2N$. 
\vspace{-20pt}
\begin{center}
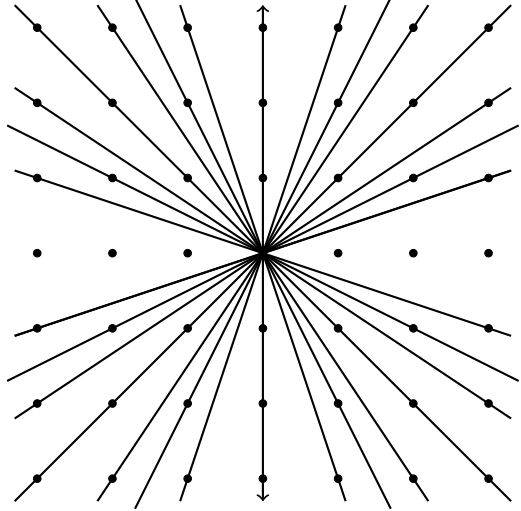
\begin{figure}[h!]
\begin{tikzpicture}[scale=1]
\draw [<->, thick] (0,-3.3) -- (0,3.3);
\filldraw (0,0) circle (0.05cm);
\filldraw (1,0) circle (0.05cm);
\filldraw (2,0) circle (0.05cm);
\filldraw (3,0) circle (0.05cm);
\filldraw (-1,0) circle (0.05cm);
\filldraw (-2,0) circle (0.05cm);
\filldraw (-3,0) circle (0.05cm);
\filldraw (0,1) circle (0.05cm);
\filldraw (1,1) circle (0.05cm);
\filldraw (2,1) circle (0.05cm);
\filldraw (3,1) circle (0.05cm);
\filldraw (-1,1) circle (0.05cm);
\filldraw (-2,1) circle (0.05cm);
\filldraw (-3,1) circle (0.05cm);
\filldraw (0,2) circle (0.05cm);
\filldraw (1,2) circle (0.05cm);
\filldraw (2,2) circle (0.05cm);
\filldraw (3,2) circle (0.05cm);
\filldraw (-1,2) circle (0.05cm);
\filldraw (-2,2) circle (0.05cm);
\filldraw (-3,2) circle (0.05cm);
\filldraw (0,3) circle (0.05cm);
\filldraw (1,3) circle (0.05cm);
\filldraw (2,3) circle (0.05cm);
\filldraw (3,3) circle (0.05cm);
\filldraw (-1,3) circle (0.05cm);
\filldraw (-2,3) circle (0.05cm);
\filldraw (-3,3) circle (0.05cm);
\filldraw (0,-1) circle (0.05cm);
\filldraw (1,-1) circle (0.05cm);
\filldraw (2,-1) circle (0.05cm);
\filldraw (3,-1) circle (0.05cm);
\filldraw (-1,-1) circle (0.05cm);
\filldraw (-2,-1) circle (0.05cm);
\filldraw (-3,-1) circle (0.05cm);
\filldraw (0,-2) circle (0.05cm);
\filldraw (1,-2) circle (0.05cm);
\filldraw (2,-2) circle (0.05cm);
\filldraw (3,-2) circle (0.05cm);
\filldraw (-1,-2) circle (0.05cm);
\filldraw (-2,-2) circle (0.05cm);
\filldraw (-3,-2) circle (0.05cm);
\filldraw (0,-3) circle (0.05cm);
\filldraw (1,-3) circle (0.05cm);
\filldraw (2,-3) circle (0.05cm);
\filldraw (3,-3) circle (0.05cm);
\filldraw (-1,-3) circle (0.05cm);
\filldraw (-2,-3) circle (0.05cm);
\filldraw (-3,-3) circle (0.05cm);
\draw [thick] (-3*1.1, -1*1.1) -- (3*1.1,1*1.1);
\draw [thick] (-3*1.1, 1*1.1) -- (3*1.1,-1*1.1);
\draw [thick] (-3*1.1, -1*1.1) -- (3*1.1,1*1.1);
\draw [thick] (-2*1.7, -1*1.7) -- (2*1.7,1*1.7);
\draw [thick] (-2*1.7, 1*1.7) -- (2*1.7,-1*1.7);
\draw [thick] (-3*1.1, -3*1.1) -- (3*1.1,3*1.1);
\draw [thick] (-3*1.1, 3*1.1) -- (3*1.1,-3*1.1);
\draw [thick] (-3*1.1, -2*1.1) -- (3*1.1,2*1.1);
\draw [thick] (-3*1.1, 2*1.1) -- (3*1.1,-2*1.1);
\draw [thick] (-2*1.1, -3*1.1) -- (2*1.1,3*1.1);
\draw [thick] (-2*1.1, 3*1.1) -- (2*1.1,-3*1.1);
\draw [thick] (-1*1.7, -2*1.7) -- (1*1.7,2*1.7);
\draw [thick] (-1*1.7, 2*1.7) -- (1*1.7,-2*1.7);
\draw [thick] (-1*1.1, -3*1.1) -- (1*1.1,3*1.1);
\draw [thick] (-1*1.1, 3*1.1) -- (1*1.1,-3*1.1);
\end{tikzpicture}
\caption{Covering lattice points outside the $x-$axis with lines.}
\end{figure}
\end{center}
\vspace{-25pt}
There are $\sim N^2$ such lines (no improvement over the trivial bound is possible because asymptotically $6/\pi^2$ of all lattice points have coprime coordinates). 
Altogether, this implies
$$   \sum_{ \|k\| \leq N}{ |\widehat{f}(k)|^2} -   \sum_{|k_1| \leq N}{ | \widehat{f}(k_1,0) |^2}  \lesssim N^2 \max_{|\gamma| \leq 2N}{~~ \left|  \frac{1}{|\gamma|} \int_{\gamma}{f d\mathcal{H}^1}  \right|^2}.$$
Using our assumption on the amount of $L^2-$mass on frequencies $\mathbb{Z} \times \left\{0\right\}$, we can derive the weaker inequality
$$   \sum_{ \|k\| \leq N}{ |\widehat{f}(k)|^2} -  \frac{\|f\|_{L^2}^2}{2}    \lesssim N^2 \max_{|\gamma| \leq 2N}{~~ \left|  \frac{1}{|\gamma|} \int_{\gamma}{f d\mathcal{H}^1}  \right|^2}.$$

We now bound the left-hand side in a way that implies a natural choice for $N$ by showing
$$ \sum_{\|k\| \geq \|\nabla f\|_{L^2}\| f\|^{-1}_{L^2}}{|\widehat{f}(k)|^2} \leq \frac{\|f\|^2_{L^2}}{4}.$$
If the inequality was false, then we could reach a contradiction by using
\begin{align*}
  \| \nabla f\|^2_{L^2} &= 4\pi^2 \sum_{k \in \mathbb{Z}^2}{|k|^2|\widehat{f}(k)|^2} \geq   4\pi^2 \sum_{\|k\| \geq \|\nabla f\|_{L^2}\| f\|^{-1}_{L^2} }{|k|^2|\widehat{f}(k)|^2} \\
&\geq  4\pi^2\frac{\|\nabla f\|^2_{L^2(\mathbb{T}^2)^2}}{\| f\|^2_{L^2(\mathbb{T}^2)^2}}\sum_{\|k\| \geq \|\nabla f\|_{L^2}\| f\|^{-1}_{L^2}}{|\widehat{f}(k)|^2} 
\geq \pi^2  \| \nabla f\|^2_{L^2(\mathbb{T}^d)}.
\end{align*}
We thus derive, for the choice
$$ N =  \frac{\|\nabla f\|_{L^2}}{\| f\|^{}_{L^2}}, \qquad \mbox{that}  \qquad \sum_{ \|k\| \leq N}{ |\widehat{f}(k)|^2} -  \frac{\|f\|_{L^2}^2}{2}  \geq \frac{\|f\|_{L^2}^2}{4}.$$

Altogether, we have thus shown the existence of a relatively short geodesic satisfying
$$ \max_{|\gamma| \leq \|\nabla f\|_{L^2}\| f\|^{-1}_{L^2} }{~~ \left|  \frac{1}{|\gamma|} \int_{\gamma}{f~~ d\mathcal{H}^1}  \right|} \gtrsim \frac{\left\| f \right\|^2_{L^2}}{ \left\| \nabla f \right\|_{L^2}}.$$
This concludes the first part of the proof.\\

 The second part of the proof consists of showing that sufficiently large geodesics yield integrals that are always smaller than the lower bound we just obtained. The combination
of those two facts then yields the desired outcome.
Consider the 'long' closed geodesic $\gamma(t) = (at, bt + c)$, long meaning that $\sqrt{a^2+b^2}$ is large, as embedded in a one-parameter family 
$\gamma_c$ indexed by $c$ (we can assume that $a \neq 0 \neq b$ since the geodesics has length 1 otherwise). We will assume w.l.o.g. that $a \geq b$: if that is not the case, then we run the subsequent argument on $\tilde f (x,y) = f(y,x)$ (the subsequent argument is only using bounds on the Fourier coefficients of $f$ that are radial and has no preferred directions).
We know that, averaged over $c$, the mean value of the average value over $\gamma_c$ is 0 and we want to show that it can never be very large. Moreover, since
$$ \gamma(t) = (at, bt + c),$$
there is a periodicity in $c$ with period $\lesssim |\gamma|^{-1}$ because the geodesic crosses the $x-$axis at least $\sim |\gamma|$ times at equally spaced intervals because $a \geq b$.  We now use a standard estimate for periodic functions $g:\mathbb{T} \rightarrow \mathbb{R}$ with mean value 0: there exists at least one point for which $g(x_0) = 0$. Then, however
\begin{align*}
 g(x)^2 = g(x)^2  - g(x_0)^2 = \int_{x_0}^{x}{\frac{d}{dy} g(y)^2 dy} =  2  \int_{x_0}^{x}{ g(y) g'(y) dy} \leq 2\|g\|_{L^2(x_0, x)}\|g'\|_{L^2(x_0, x)}.
\end{align*}
The periodicity in $c$ implies that we can pick the point at distance $\lesssim |\gamma|^{-1}$ from $x$ and thus
$$ 2\|g\|_{L^2(x_0, x)}\|g'\|_{L^2(x_0, x)} \lesssim \frac{1}{|\gamma|} \|g\|_{L^2(\mathbb{T})}\|g'\|_{L^2(\mathbb{T})}.$$
Therefore, since $x$ was completely arbitrary,
$$ \|g\|_{L^{\infty}(\mathbb{T})} \lesssim \frac{1}{|\gamma|^{1/2}}\|g\|^{1/2}_{L^{2}(\mathbb{T})} \|g'\|^{1/2}_{L^{2}(\mathbb{T})}.$$
This implies that 
\begin{align*}
\max_{0 \leq c \leq 1}{ \frac{1}{|\gamma_{a,b,c}|} \left| \int_{\gamma_{a,b,c}}{f~~d\mathcal{H}^1}  \right|} &= \left\|  \sum_{ k_1 a + k_2 b = 0}{ \widehat{f}(k) e^{2\pi i k_2 c} } \right\|_{L^{\infty}[0,1]} \\
&\lesssim   \frac{1}{|\gamma|^{1/2}} \left(  \sum_{ k_1 a + k_2 b = 0}{| \widehat{f}(k_1, k_2)|^2 } \right)^{1/4}
  \left(  \sum_{ k_1 a + k_2 b = 0}{\left\|k_2\right\|^2 | \widehat{f}(k_1,k_2)|^2 } \right)^{1/4}.
\end{align*}
We bound these quantities by invoking uniform bounds on the Fourier coefficients of $C^s$ functions: for any $s \geq 1$
$$ | \widehat{f}(k)| \lesssim_{s} \frac{ \max_{|\alpha| = s}{\left\|  \partial_{\alpha} f  \right\|_{L^{1}}}}{ \|k\|^s}.$$
Using this uniform estimate, we obtain the following bound for the first term
\begin{align*}
 \sum_{ k_1 a + k_2 b = 0}{| \widehat{f}(k_1, k_2)|^2} &= \sum_{ d \in \mathbb{Z} \atop d \neq 0}{| \widehat{f}(- d b, d a)|^2} \\
&\lesssim_{s} \left( \max_{|\alpha| = s}{\left\| \partial_{\alpha} f  \right\|_{L^{1}}}\right)^2 \sum_{d \in \mathbb{Z} \atop d \neq 0}{ \frac{1}{d^{2s}} \frac{1}{(a^2+b^2)^{s}}}\\
&\lesssim_s \frac{ \left( \max_{|\alpha| = s}{\left\| \partial_{\alpha} f  \right\|_{L^{1}}}\right)^2}{(a^2+b^2)^{s}}.
\end{align*}
The same kind of estimate can be applied to the second term and results in
\begin{align*}
 \sum_{ k_1 a + k_2 b = 0}{ |k_2|^2| \widehat{f}(k_1, k_2)|^2} &= \sum_{ d \in \mathbb{Z} \atop d \neq 0}{d^2 a^2 | \widehat{f}(- d b, d a)|^2} \\
&\lesssim_{s} \left( \max_{|\alpha| = s}{\left\| \partial_{\alpha} f  \right\|_{L^{1}}}\right)^2 \sum_{d \in \mathbb{Z} \atop d \neq 0}{ \frac{1}{d^{2(s-1)}} \frac{1}{(a^2+b^2)^{s-1}}}\\
&\lesssim_s \frac{ \left( \max_{|\alpha| = s}{\left\| \partial_{\alpha} f  \right\|_{L^{1}}}\right)^2}{(a^2+b^2)^{s-1}},
\end{align*}
where the last step requires $s \geq 2$.
For $s \geq 2$, these two bounds imply
$$ \left\|  \sum_{ k_1 a + k_2 b = 0}{ \widehat{f}(k) e^{2\pi i k_2 c} } \right\|_{L^{\infty}[0,1]}  \lesssim_s  \frac{ \max_{|\alpha| = s}{\left\| \partial_{\alpha} f  \right\|_{L^{1}}}}{(a^2+b^2)^{\frac{s}{2} }}.$$
The length of a closed geodesic indexed by $\gamma(t) = (at, bt + c)$ is, since $a,b$ have no common divisors, given by $\sqrt{a^2+b^2}$. Altogether, this implies the uniform estimate
$$  \left|  \frac{1}{|\gamma|} \int_{\gamma}{f~~ d\mathcal{H}^1}  \right| \lesssim_s  \frac{ \max_{|\alpha| = s}{\left\| \partial_{\alpha} f  \right\|_{L^{1}}}}{|\gamma|^{s}}$$
However, this allows for finding a critical length beyond which geodesics are bound to be suboptimal by solving for
$$ \frac{ \max_{|\alpha| = s}{\left\| \partial_{\alpha} f  \right\|_{L^{1}}}}{|\gamma|^{s }} \lesssim_s \frac{\left\| f \right\|^2_{L^2}}{ \left\| \nabla f \right\|_{L^2}}$$
which yields the desired result.
\end{proof}

\subsection{Proof of the Corollary}
\begin{proof} The main result states that
$$ |\gamma|^{s} \lesssim_s  \left( \max_{|\alpha| = s}{ \| \partial_{\alpha} f \|_{L^{1}(\mathbb{T}^2)}} \right) \| \nabla f \|_{L^2}^{} \|f\|_{L^2}^{-2}.$$
For functions 
$$ f(x) = \sum_{\|k\| \leq N}{\widehat{f}(k) e^{2\pi i \left\langle k, x \right\rangle}},~ \mbox{the quantity} \qquad \max_{|\alpha| = s}{\left\|  \partial_{\alpha} f  \right\|_{L^{1}}}$$
can be controlled fairly well: let us assume that $\alpha =(\alpha_1, \alpha_2)$ with $\alpha_1 + \alpha_2 = s$. Then, using the triangle inequality and the Cauchy-Schwarz inequality,
\begin{align*}
  \left\|  \partial_{\alpha} f  \right\|_{L^{1}}  = \left\|  \sum_{\|k\| \leq N}{ k_1^{\alpha_1} k_2^{\alpha_2} \widehat{f}(k) e^{2\pi i \left\langle k, x \right\rangle}}  \right\|_{L^{1}}
&\leq  \sum_{\|k\| \leq N}{  |k_1|^{\alpha_1} |k_2|^{\alpha_2} |\widehat{f}(k)| } \\
&\leq N^s  \sum_{\|k\| \leq N}{ |\widehat{f}(k)| } \lesssim N^{s+1} \|f\|_{L^2}.
\end{align*}
It remains to compute the implicit constant in $\lesssim_s$ which is known \cite[Theorem 3.2.9]{grafakos} to grow at most like $c^s$ for a universal constant $c>0$. Letting $s \rightarrow \infty$
then implies the result. In fact, a stronger statement is true (and would yield a simple alternative proof): the main algebraic identity shows that all longer closed geodesics  are actually orthogonal to $f$ and always yield 0.
\end{proof}

\end{document}